\documentclass[11pt,reqno]{amsart}
\usepackage[margin=1in,letterpaper]{geometry}
\usepackage{graphicx}
\usepackage{amssymb}
\usepackage{amsthm}
\usepackage{mathrsfs}
\usepackage{stmaryrd}
\usepackage{accents} 
\usepackage{enumitem} 
\usepackage{subcaption}

\makeatletter\let\over\@@over\makeatother

\numberwithin{equation}{section}
\theoremstyle{plain} 
\newtheorem{theorem}{Theorem}[section] 
\newtheorem{proposition}[theorem]{Proposition}

\theoremstyle{remark}

\theoremstyle{definition}

\newcommand{\be}{\begin{equation}}
\newcommand{\ee}{\end{equation}}%
\newcommand{\bse}{\begin{subequations}}
\newcommand{\ese}{\end{subequations}}

\newcommand{\jump}[1]{\left\llbracket{#1}\right\rrbracket}

\newcommand{\kernel}{\operatorname{ker}}
\newcommand{\linspan}{\operatorname{span}}

\newcommand{\R}{\mathbb{R}} 
\newcommand{\placeholder}{\,\cdot\,}

\newcommand{\n}[2][]{#1\lVert #2 #1\rVert}
\newcommand{\abs}[2][]{#1\lvert #2 #1\rvert}

\newcommand{\ona}{\textup{~on~}}

\newcommand{\bdd}{\mathrm{b}}       
\newcommand{\loc}{{\mathrm{loc}} }     

\newcommand\F{\mathscr F}    

\newcommand{\cm}{{\mathscr C}}  
\newcommand{\iftcm}{{\mathscr C}}

\newcommand\prineigenvaluepm{\sigma_{0}^\pm}
\newcommand\prineigenvalue{\sigma_{0}}

\newcommand\fluidD{\mathscr{D}}
\newcommand\fluidS{\mathscr{S}}

\newcommand\Lip{\operatorname{Lip}}
\newcommand\lambdastar{\lambda_+} 

\title{Large-amplitude internal fronts in two-fluid systems}

\author[R. M. Chen]{Robin Ming Chen}
\address{Department of Mathematics, University of Pittsburgh, Pittsburgh, PA 15260} 
\email{mingchen@pitt.edu}  

\author[S. Walsh]{Samuel Walsh}
\address{Department of Mathematics, University of Missouri, Columbia, MO 65211} 
\email{walshsa@missouri.edu}

\author[M. H. Wheeler]{Miles H. Wheeler}
\address{Department of Mathematical Sciences, University of Bath, Bath BA2 7AY, United Kingdom}
\email{mw2319@bath.ac.uk}

\begin{document}

\begin{abstract} 
In this announcement, we report results on the existence of families of large-amplitude internal hydrodynamic bores.  These are traveling front solutions of the full two-phase incompressible Euler equation in two dimensions.  The fluids are bounded above and below by flat horizontal walls and acted upon by gravity.  We obtain continuous curves of solutions to this system that bifurcate from the trivial solution where the interface is flat.     Following these families to the their extreme,  the internal interface either overturns, comes into contact with the upper wall, or develops a highly degenerate ``double stagnation'' point.

Our construction is made possible by a new abstract machinery for global continuation of monotone front-type solutions to elliptic equations posed on infinite cylinders.  This theory is quite robust and, in particular, can treat fully nonlinear equations as well as quasilinear problems with transmission boundary conditions.  \end{abstract}

\maketitle

\section{Introduction} \label{intro bores section}

The world's oceans are stratified in the sense that the fluid density increases with depth. While small in relative terms, this density variation can dramatically affect the dynamics and, in particular, allows for the formation of large scale \emph{internal waves} that remain coherent over long distances. In many settings there are two regions with nearly constant density separated by a thin layer, called the \emph{pycnocline}, where density gradients are large. This permits the system to be modeled as two constant density fluids with different densities, divided by a sharp interface along which  waves can propagate. Unlike surface waves in a homogeneous density fluid, these internal waves can take the form of fronts or (smooth) \emph{hydrodynamical bores}. These are steady solutions where the internal interface is asymptotically flat both upstream and downstream of the wave but with different heights.

Let us restrict attention to the simplest configuration where the two fluid layers are irrotational and bounded from above and below by rigid flat boundaries as shown in Figure~\ref{basic bore figure}. There is an extensive applied literature on this problem, mostly centered around linear or weakly nonlinear model equations which are valid only for small amplitudes \cite{helfrich2006review}, as well as a growing body of rigorous results. For bores in the full nonlinear equations, the first rigorous existence results date back to the work of Amick and Turner~\cite{amick1989small}, confirming formal predictions based on the weakly nonlinear extended Korteweg--de Vries equation. Alternative proofs have subsequently been given using different methods by Mielke~\cite{mielke1995homoclinic},  Makarenko~\cite{makarenko1992bore}, and the authors~\cite{chen2019center}. 
 
In this announcement, we report the first construction of genuinely large-amplitude bores. One can no longer expect to base such an analysis on a well-chosen model equation, and instead we rely on a new abstract global bifurcation theory tailored to front-type solutions of elliptic equations in cylindrical domains \cite{chen2020global}. 

\subsection{Formulation and existence theory}

\begin{figure}[tb]
  \centering
  \includegraphics[scale=0.85]{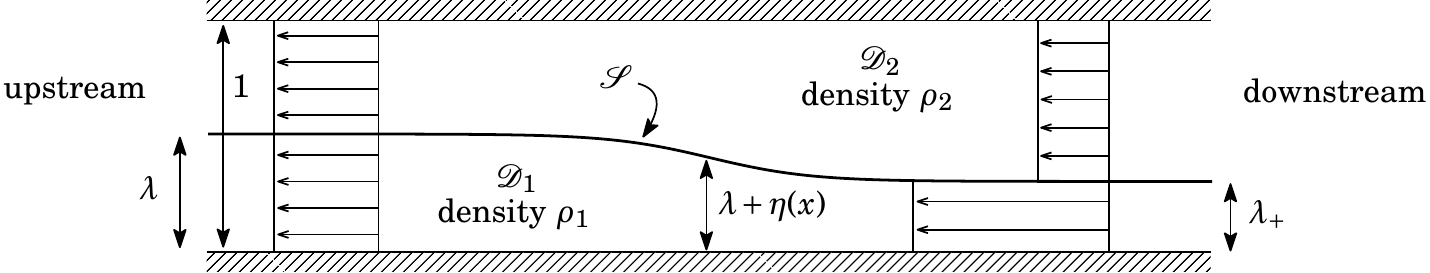}
  \caption{A bore in a two layer fluid.}
  \label{basic bore figure}
\end{figure}

The problem can be mathematically formulated as follows. The unknown interface $\fluidS = \{ (x,y) : y=\eta(x) \}$ separates two open fluid regions $\fluidD_1$ and $\fluidD_2$ as shown in Figure~\ref{basic bore figure}. Here the lower region $\fluidD_1$ has constant density $\rho_1 > 0$, and is bounded below by a rigid barrier at height $y=-\lambda$. The upper region $\fluidD_2$ is likewise bounded above by a rigid barrier at $y=1-\lambda$ and has constant density $0 < \rho_2 < \rho_1$.  Note that the total height of the channel is normalized to $1$.

\begin{subequations}\label{eqn:stream}
 Assuming incompressibility, the velocity field in each fluid is given by $(\partial_y \psi, -\partial_x \psi)$ for some \emph{stream function} $\psi$ satisfying
  \begin{equation} \label{eqn:psi harmonic}
    \Delta \psi = 0 \quad \text{in } \fluidD_1 \cup \fluidD_2
  \end{equation}
  together with the so-called kinematic boundary conditions
  \begin{align}
    \label{eqn:stream:kin}
    \begin{alignedat}{2}
      \psi &= 0 &\qquad& \text{ on } \fluidS,\\
      \psi &= \lambda &\quad& \text{ on } y=-\lambda,  \\
      \psi &= \lambda-1 &\quad& \text{ on } y=1-\lambda, 
    \end{alignedat}
  \end{align}
  and the dynamic boundary condition
  \begin{equation}
    \label{eqn:stream:dynamic}
    \frac 12 \jump{\rho\abs{\nabla\psi}^2}
    + \frac {\jump\rho}{F^2} y = \frac {\jump\rho}2  \qquad  \ona \fluidS,
  \end{equation}
  where $F>0$ is a dimensionless parameter called the \emph{Froude number} and
  $\jump{\placeholder}=(\placeholder)_2-(\placeholder)_1$ denotes the jump of a quantity across the interface $\fluidS$. These are supplemented with the asymptotic conditions
  \begin{align}
    \label{eqn:stream:asym}
    \nabla\psi \to (0,-1),
    \quad 
    \eta \to 0
    \qquad &\text{ as } x \to -\infty,\\
    \eta \to \lambdastar - \lambda \ne 0 \qquad \quad &\text{ as } x \to +\infty
  \end{align}
  as shown in Figure~\ref{basic bore figure}.
\end{subequations}

We are interested in classical solutions to \eqref{eqn:stream} which enjoy the regularity
\begin{equation*}
  \psi \in C_\bdd^{2+\alpha}(\overline{\fluidD_1}) \cap C_\bdd^{2+\alpha}(\overline{\fluidD_2}) \cap C_\bdd^0(\overline{\fluidD_1 \cup \fluidD_2}), \qquad \eta \in C_\bdd^{2+\alpha}(\mathbb{R}),
\end{equation*}
for a fixed $\alpha \in (0,1)$, and where the subscript `b' indicates uniform boundedness.  As is well known in the literature on internal waves, this situation is only possible provided the Froude number $F$ and constant $\lambdastar$ are given explicitly by 
\begin{align}
  \label{conjugate flow conclusion}
  F^2 = \frac{\sqrt{\rho_1}-\sqrt{\rho_2}}{\sqrt{\rho_1}+\sqrt{\rho_2}},
  \qquad 
  \lambdastar = \frac {\sqrt{\rho_1}}{\sqrt{\rho_1}+\sqrt{\rho_2}},
\end{align}
see for instance \cite[Appendix~A]{Laget1997interfacial}. Note that we are working in units where the height of the channel is the length scale, and the upstream (relative) velocity is the velocity scale. 

Our first theorem is the following global bifurcation result \cite{chen2020global}.

\begin{theorem}[Large-amplitude bores] \label{global bore theorem} 
  Fix $\alpha \in (0,1)$ and densities $0 < \rho_2 < \rho_1$.  There exist $C^0$ curves
  \begin{equation*}
    \cm^\pm = \left\{ (\psi(s), \eta(s), \lambda(s)) : \pm s \in (0,\infty) \right\}
  \end{equation*}
  of classical solutions to the internal wave problem \eqref{eqn:stream}--\eqref{conjugate flow conclusion} with the following properties.
  \begin{enumerate}[label=\rm(\alph*)]
  \item \label{global bore monotone part} \textup{(Strict monotonicity)} Each solution on $\cm^\pm$ is a strictly monotone bore:
    \begin{equation*}
      \begin{alignedat}{2}
        \pm \partial_x \eta(s) &< 0 & \quad &  \textup{on } \mathbb{R}, \\
        \pm \partial_x \psi(s) &> 0 & \quad & \textup{in } \fluidD_1(s) \cup \fluidD_2(s) \cup \fluidS(s), \\
        \partial_y \psi(s) & < 0 & \quad & \textup{in } \overline{\fluidD_1(s) \cup \fluidD_2(s)}.
      \end{alignedat} 
    \end{equation*}
  \item \label{global bore limit part} \textup{(Stagnation limit)} Following $\cm^\pm$, we encounter waves that are arbitrarily close to having a horizontal stagnation point on the internal interface:
    \begin{equation}
      \lim_{s \to \pm\infty} \sup_{\fluidS(s)} \partial_y\psi_i(s) = 0, \qquad \textup{for $i = 1$ or $2$}.  \label{stagnation limit} 
    \end{equation}
  \item \label{global bore laminar part} \textup{(Laminar origin)} Both $\cm^-$ and $\cm^+$ emanate from the same laminar solution in that 
    \begin{equation*}
      \eta(s) \to 0, \quad \nabla \psi(s) \to (0,-1), \quad \lambda(s) \to {\lambdastar\pm} \qquad \textup{as } s \to 0\pm.
    \end{equation*}
  \end{enumerate}
\end{theorem}

An overview of the proof of Theorem~\ref{global bore theorem} is given at the end of Section~\ref{abstract results section}.

\begin{figure}[tb]
  \centering
  \includegraphics[scale=0.85]{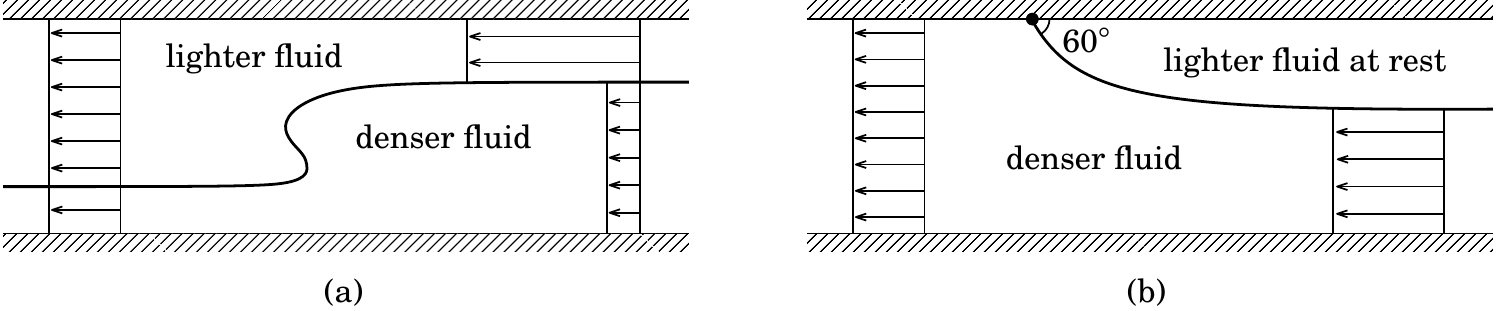}
  \caption{Limiting configurations of monotone bores.
   }
  \label{bores alternatives figure}
\end{figure}

\subsection{Overhanging water waves}
The next result characterizes the limiting form of the profile along $\cm^\pm$.  

\begin{theorem}[Limiting interface] \label{limit eta theorem} 
  \hfill
  \begin{enumerate}[label=\rm(\alph*)]
  \item \textup{(Overturning or singularity)} \label{elevation part} In the limit along $\cm^-$, either the interface \emph{overturns} in that
    \begin{equation*}
      \limsup_{s \to -\infty} \| \partial_x \eta(s) \|_{L^\infty(\mathbb{R})} = \infty, 
   \end{equation*}
  or it becomes \emph{singular} in that we can extract a translated subsequence 
  \[ \eta(s) \longrightarrow \eta^* \in \Lip(\mathbb{R}) \quad \textrm{in $C_\loc^\varepsilon$ for all } \varepsilon \in (0,1) \]
   such that $\{ y < \eta^*(x) \}$ simultaneously fails to satisfy both an interior sphere and exterior sphere condition at a single point on its boundary. 
  \item \textup{(Overturning or contact)} \label{depression part} Following $\cm^+$, either the interface overturns or it comes into contact with the upper wall:
  \begin{equation*}
    \limsup_{s \to \infty} \lambda(s) = 1 \qquad \textrm{or} \qquad  \limsup_{s \to \infty} \| \partial_x \eta(s) \|_{L^\infty(\mathbb{R})} = \infty.
  \end{equation*}
  \end{enumerate}
\end{theorem}
We conjecture that the singularity alternative in \ref{elevation part} can be eliminated with further analysis, and hence that overturning occurs. This would be consistent with numerical work of Dias and Vanden-Broeck~\cite{dias2003internal}, which suggests that the interfaces along $\cm^-$ overturn while those along $\cm^+$ come into contact with the upper wall as shown in Figure~\ref{bores alternatives figure} to form a so-called ``gravity current''.

The conclusion of Theorem~\ref{limit eta theorem} should be compared with earlier work on overturning or \emph{overhanging} water waves. These are waves for which the interface $\fluidS$ ceases to be a graph. Since we are speaking of traveling wave solutions, this configuration persists for all time, a quite strange phenomenon that has sparked intense mathematical interest.

There are three main physical effects to be considered: gravity, surface tension forces along the interface, and vorticity, either continuously distributed throughout the fluid or concentrated into a internal vortex sheet as in \eqref{eqn:stream}. When only gravity is present, overturning waves cannot exist \cite{spielvogel:overturning,amick1987bounds,toland1996stokes}. When only surface tension forces or capillarity is present, on the other hand, there is an exact family of overhanging periodic solutions due to Crapper \cite{crapper1957exact}. These waves have subsequently been perturbed by Akers, Ambrose, and Wright \cite{akers2014gravity} and C\'ordoba, Enciso, and Grubic \cite{cordoba2016existence} to obtain overhanging gravity-capillary waves where the (dimensionless) gravity is small. Very recently, another family of explicit overhanging waves has been discovered \cite{dyachenko2019folds,hur2020crapper,hur2020exact} where only vorticity is present. 

In the above existence results, gravity is either neglected or treated as a small parameter. However, numerical work has shown that overhanging waves can also exist when gravity effects are $O(1)$, and in particular in the absence of surface tension forces \cite{pullin1988finite,turner1988broadening}. Rigorous verification of these results remains an outstanding open problem, but significant progress has been made by Constantin, Varvaruca, and Strauss \cite{constantin2016global}. They used global bifurcation theory to construct a continuous curve of solutions, which in principle are allowed to overturn. Based on subsequent numerics \cite{dyachenko2019folds,dyachenko2019stokes} it is conjectured that these curves indeed contain overhanging waves, but a rigorous proof has so far been elusive. We mention related work on periodic internal gravity waves \cite{sun2001existence,maklakov2020note} and internal gravity-capillary waves \cite{ambrose2016global}. 

Theorem~\ref{limit eta theorem}\ref{elevation part} is tantalizingly close to a proof of overturning. The only other possibility is a highly degenerate type of singularity that has not been observed in numerics \cite{dias2003internal}. Moreover, there is hope that such singularities could be ruled out through a completely local analysis. By comparison, the global bifurcation results in \cite{constantin2016global,sun2001existence,ambrose2016global} allow for a wide range of possibilities. The price we pay for this apparent advantage is twofold. First, we work with a reformulation of \eqref{eqn:stream} that degenerates as an overturning wave is approached. This allows us to detect overturning more easily, but prevents us from continuing further to obtain truly overhanging waves. Second, we construct bores rather than periodic or solitary waves. This introduces serious difficulties related to the unboundedness of the fluid domain as well as the lack of symmetry for the solutions. Once these considerable obstacles have been overcome, however, we find ourselves with more concrete information about the solutions than would be available in the periodic or solitary wave cases.

\section{Global bifurcation of monotone fronts}

\subsection{Motivation from second-order ODEs}

As mentioned above, Theorem~\ref{global bore theorem} is obtained through a much more general set of results on the global bifurcation of monotone fronts in elliptic PDE.  Before presenting those ideas, let us briefly discuss the setting of second-order ODEs where it is easier to construct concrete examples.  

In fact, an equation of this type is frequently used as a simplified model for the internal wave system \eqref{eqn:stream}.  Under the assumption that the waves are long (in some appropriate sense) but not necessarily small amplitude, Miyata \cite{miyata1988long} and Choi--Camassa \cite{choi1996long,choi1999fully} independently derived a time-dependent PDE related to the Serre--Green--Naghdi system. Referred to as the MCC equation, this model reduces to the extended Korteweg--de Vries equation mentioned in the introduction in the small-amplitude limit, but it is far more accurate for waves of moderate and even large amplitude. With our current notation, the MCC equation reads
\begin{equation}
  \zeta_x^2 = \frac{3 \zeta^2}{2F^2} \frac{(\lambda+\zeta)(1-\lambda-\zeta+F^2) \rho_2-(1-\lambda-\zeta)(\lambda+\zeta-F^2) \rho_1}{(1-\lambda)^2(\lambda+\zeta)\rho_2+\lambda^2(1-\lambda-\zeta)\rho_1}  
  \label{mcc equation} 
\end{equation}
in integrated form, where here we write $\zeta$ rather than $\eta$ for the deflection of the interface to emphasize the distinction with the full system \eqref{eqn:stream}. In differentiated form, \eqref{mcc equation} can be written as
\begin{equation}
  \ddot{\zeta} + V_z(\zeta,\lambda) = 0  \label{potential equation} 
\end{equation}
for an explicit $V = V(z,\lambda)$ that is analytic in its arguments and where dot denotes derivative in $x$. Here we are viewing the densities $\rho_1,\rho_2$ as fixed and the Froude number $F$ as given by \eqref{conjugate flow conclusion}, so that the upstream depth $\lambda$ of the lower fluid layer is the only parameter. A bore now corresponds to a heteroclinic orbit of \eqref{potential equation} connecting two distinct equilibria.

For general equations of the form \eqref{potential equation} and a given pair of equilibria, it is relatively straightforward to formulate general conditions which guarantee the existence of heteroclinic orbits. For instance we have the following.

\begin{proposition}\label{ode proposition} Consider the second-order ODE \eqref{potential equation}.  Suppose that for a fixed parameter $\lambda_0$, there are two distinct rest points $Z_-(\lambda_0)$ and $Z_+(\lambda_0)$ that are \emph{conjugate} in that
  \begin{equation}
    V(Z_-(\lambda_0),\lambda_0) = V(Z_+(\lambda_0),\lambda_0). \label{ode conjugate} 
  \end{equation}
Assume also that the potential satisfies a heteroclinic nondegeneracy condition
  \begin{equation}
    V(z, \lambda_0) < V(Z_\pm(\lambda_0),\lambda_0) \qquad \textrm{for $z$ between $Z_+(\lambda_0)$ and $Z_-(\lambda_0)$},   \label{ode hetero nondegen condition} 
  \end{equation}
and spectral nondegeneracy condition 
  \begin{equation}
    V_{zz}(Z_-(\lambda_0), \lambda_0),~V_{zz}(Z_+(\lambda_0),\lambda_0) < 0. \label{ode spectral nondegen condition}
  \end{equation} 
Then there exists a solution $(\zeta_0, \lambda_0)$ to \eqref{potential equation} with $\zeta_0(x) \to Z_\pm(\lambda_0)$ as $x \to \pm\infty$.
\end{proposition}

Stated simply, the problem of finding heteroclinic solutions to the ODE \eqref{potential equation} amounts to verifying the existence of conjugate rest points of $V(\placeholder, \lambda_0)$ satisfying a type of heteroclinic nondegeneracy condition \eqref{ode hetero nondegen condition} and a spectral nondegeneracy condition \eqref{ode spectral nondegen condition}.  
We can moreover consider the case when there is a smooth family of conjugate rest points $Z_+(\lambda)$ and $Z_-(\lambda)$ that satisfy \eqref{ode conjugate}--\eqref{ode spectral nondegen condition} for $\lambda$ in a neighborhood of $\lambda_0$.  It is not hard to see that there will then exist a local curve $\mathscr{K}_\loc$ of heteroclinic orbits bifurcating from $(\zeta_0, \lambda_0)$.  Clearly, one can continue this curve at least as far as the above hypotheses are satisfied along it.  

Applying Proposition~\ref{ode proposition} to the MCC model \eqref{mcc equation}, we find that for any $\lambda \in (0,1)$ there is always a unique smooth heteroclinic orbit connecting the rest points $Z_-(\lambda) = 0$ and $Z_+(\lambda)=\lambda_+$, where here $\lambda_+$ is given by \eqref{conjugate flow conclusion}. Recall that for the full problem, numerical evidence \cite{dias2003internal} suggests that some bores are instead overturning. Such waves would violate the long-wave assumption made in the derivation of \eqref{mcc equation}, and so this discrepancy is to be expected.

\subsection{Monotone fronts solutions to elliptic PDE}

Keeping in mind the above discussion, consider now the following (fully) nonlinear PDE:
\begin{equation}
 \label{fully nonlinear elliptic pde}
  \left\{ \begin{aligned}
    {A}(y, u, \nabla u, D^2 u,  \lambda)  & = 0  \qquad \textrm{in } \Omega, \\
    B(y, u, \nabla u, \lambda) & = 0 \qquad  \textrm{on } \Gamma_1, \\
    u & = 0 \qquad \textrm{on } \Gamma_0,
  \end{aligned} \right.
\end{equation}
where $\lambda \in \R$ is a parameter, and the domain $\Omega = \mathbb{R} \times \Omega^\prime$ is an infinite cylinder with bounded base $\Omega^\prime \subset \mathbb{R}^{d-1}$.  For simplicity, assume that $\Omega$ is connected with a $C^{2+\alpha}$ boundary $\partial \Omega = \Gamma_0 \cup \Gamma_1$, for a fixed $\alpha \in (0,1)$ and such that $\Gamma_0 \cap \Gamma_1 = \emptyset$.   
Points in $\Omega$ will be denoted $(x,y)$, where $x \in \mathbb{R}$ and $y \in \Omega^\prime$.

We assume that ${A}$ and ${B}$ are real analytic in all of their arguments and that \eqref{fully nonlinear elliptic pde} is uniformly elliptic with a uniformly oblique boundary condition on $\Gamma_1$.  Through the Dubreil-Jacotin transform, the internal waves problem \eqref{eqn:stream} can be rewritten roughly in this form with upstream layer depth ratio as the parameter.  In fact, the dynamic condition \eqref{eqn:stream:dynamic} will lead to a nonlinear transmission problem, but this can be handled through a small modification.  
  
Define a \emph{front} to be a solution $(u, \lambda)$ of \eqref{fully nonlinear elliptic pde} that enjoys the classical regularity $u \in C^{2+\alpha}_\bdd(\overline{\Omega})$ and has distinct point-wise limits as $x \to -\infty$ and $x \to +\infty$; thinking of water waves, we call these the upstream and downstream states, respectively.  From the structure of the equation, one can prove that they are in fact $x$-independent solutions of \eqref{fully nonlinear elliptic pde}.  
We call a front \emph{monotone} provided $\partial_x u \leq 0$ (or $\partial_x u \geq 0$) in $\Omega$, and \emph{strictly monotone} if $\partial_x u < 0$ (or $\partial_x u > 0$) in $\Omega \cup \Gamma_1$.  

Fronts are the PDE analogues of heteroclinic solutions to the ODE \eqref{potential equation} with the (unbounded) axial direction identified with the evolution variable.  We may then ask:  (i) under what conditions does \eqref{fully nonlinear elliptic pde} support (monotone) fronts, and (ii) do these fronts persist for non-perturbative parameter values.  
The first of these questions has been pursued by many authors.  The most common approaches include monotonicity methods \cite{berestycki1992travelling,volpert1994book} and center manifold reduction \cite{kirchgassner1982wavesolutions,mielke1988reduction}, which has been applied to our system \eqref{eqn:stream} in \cite{amick1989small,mielke1995homoclinic,chen2019center}.

Our main abstract result addresses the second problem, namely the global continuation of a given curve $\cm_\loc$ of perturbative strictly monotone fronts.  In brief, it gives conditions analogous to those of Proposition~\ref{ode proposition} under which the local curve $\cm_\loc$ can be extended to a larger curve $\cm$ of strictly monotone fronts.  These hypotheses are discussed in the next subsection.  In Section~\ref{alternatives section}, we give a sharp set of alternatives that characterize the limiting behavior as one follows the resulting global curve to its extreme.  The statement of the global bifurcation theorem is found in Section~\ref{abstract results section}.

\subsection{Hypotheses} 
In what follows, we suppose that there exists a local curve $\cm_\loc$ of strictly monotone front solutions to \eqref{fully nonlinear elliptic pde}.  To simplify the notation, it is useful to write \eqref{fully nonlinear elliptic pde} as the abstract operator equation
\[ \F(u, \lambda) = 0.\]
One can easily verify that $\F$ is real analytic as a mapping $C_\bdd^{2+\alpha}(\overline{\Omega}) \times \mathbb{R}  \to C_\bdd^{\alpha}(\overline{\Omega}) \times C_\bdd^{1+\alpha}(\Gamma_1)$.   

First, we note that the system \eqref{fully nonlinear elliptic pde} is invariant under translation in $x$, and so $\partial_x u$ lies in $\kernel{\F_u(u,\lambda)}$ by an elliptic regularity argument.    For simplicity, we assume that along the local curve the kernel is exactly one dimensional:
\begin{equation}
  \ker \F_u(u,\lambda) = \linspan \{ \partial_x u \} \qquad \textrm{for all } (u,\lambda) \in \cm_\loc. \label{kernel assumption}\tag{H1} 
\end{equation}

The next hypothesis corresponds to the spectral non-degeneracy condition \eqref{ode spectral nondegen condition} in the ODE setting.  For a monotone front $(u,\lambda)$, the Fr\'echet derivative $\F_u(u,\lambda)$ is a linear elliptic operator whose coefficients have well-defined limits as $x \to \pm\infty$.  Restricting the domain to $x$-independent functions, this gives elliptic operators on $\Omega^\prime$ that we call the \emph{transversal linearized operators} at $x = +\infty$ and $x=-\infty$. One can show that these will have principal eigenvalues that we will denote by $\prineigenvaluepm(u,\lambda)$.  Recall from elliptic theory, the principal eigenvalue is real and lies strictly to the right of the rest of the spectrum.

 In analogy to the assumption \eqref{ode spectral nondegen condition} in Proposition~\ref{ode proposition}, we focus on the situation where 
\begin{equation}
  \prineigenvalue^-(u,\lambda),~\prineigenvalue^+(u,\lambda) < 0  \qquad \textrm{for all } (u,\lambda) \in \cm_\loc. \label{spectral assumption}\tag{H2} 
\end{equation}
Observe that  \eqref{spectral assumption} is equivalent to the essential spectrum of the limiting linearized operators being properly contained in left complex half-plane $\mathbb{C}_-$.

 The final hypothesis is made with an eye towards applications.  Usually, one obtains $\cm_\loc$ through a preliminary local bifurcation argument.  A common scenario on unbounded domains is that $\cm_\loc$ originates from an $x$-independent solution to \eqref{fully nonlinear elliptic pde} that is singular in the sense that the linearized operator there fails to be Fredholm.  With that in mind, suppose that $\cm_\loc$ admits the $C^0$ parameterization
\[ 
\mathcal{C}_\loc = \left\{ \left( u(\varepsilon), \lambda(\varepsilon) \right) : 0 < \varepsilon < \varepsilon_0 \right\}  \subset \F^{-1}(0),
\]
where 
\begin{equation}
  \begin{gathered}
    (u(\varepsilon), \lambda(\varepsilon)) \to (u_0,\lambda_0) \quad \text{as } \varepsilon \to {0+}, \qquad 
    \textrm{ and } \qquad \prineigenvalue^+(u_0,\lambda_0)=0 \text{ or } \prineigenvalue^-(u_0,\lambda_0) = 0.
  \end{gathered}
  \label{local singular assumption} 
  \tag{H4} 
\end{equation}
We label this condition \eqref{local singular assumption} rather than (H3) for consistency with \cite{chen2020global}.

\subsection{Alternatives} \label{alternatives section} 
 
Taking for granted that $\cm_\loc$ can be extended, the next question is what we might encounter at the extreme of the resulting global curve. To form intuition for the PDE case, let us consider in tandem the simpler task of continuing the curve $\mathscr{K}_\loc$ of heteroclinic solutions to the ODE \eqref{potential equation}.  

An obvious possibility is that the heteroclinic orbits persist for all parameter values (hence $\lambda$ is unbounded along the curve) or that arbitrarily large fronts exist (that is, $\zeta$ is unbounded in norm).  This alternative has a straightforward translation to the PDE setting: we say that a sequence of monotone fronts $\{(u_n, \lambda_n)\}$ experiences \emph{blowup} provided that
\begin{equation}
  \n{u_n}_{C^{2+\alpha}(\Omega)} + |\lambda_n | \longrightarrow \infty.\label{gen blowup alternative} \tag{A1} 
\end{equation}
Note that in applications, it is often necessary to formulate the theory for $(u,\lambda)$ lying in an open subset of $C^{2+\alpha}(\overline{\Omega}) \times \mathbb{R}$ rather than the whole space.  For example, our analysis of internal waves supposes the absence of horizontal stagnation, which corresponds to a pointwise inequality for a certain derivative of $u$.    In  that case, the definition of blowup will include the possibility that the sequence limits to the boundary of this set.

\begin{figure}
  \centering
  \includegraphics[scale=1,page=1]{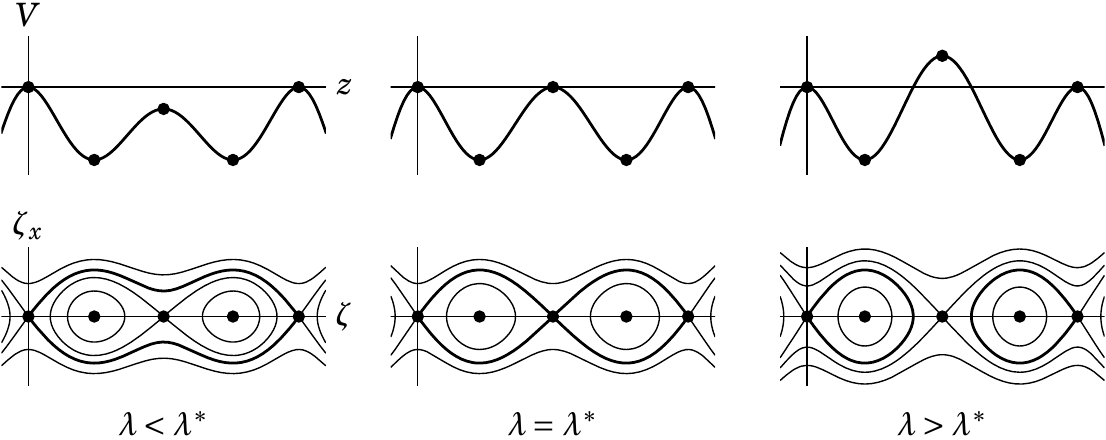}

  \caption{An ODE of the form \eqref{potential equation} that experiences a heteroclinic degeneracy \eqref{gen hetero degen}.
  } \label{hetero degen phase portrait figure}
\end{figure}

Another alternative for the ODE \eqref{potential equation} is that the heteroclinic nondegeneracy condition \eqref{ode hetero nondegen condition} is violated in the limit.  It could then happen that the heteroclinic orbit between the equilibria $Z_-(\lambda)$ and $Z_+(\lambda)$ breaks down and a new heteroclinics is born that connects one of them to an intermediate rest point as in Figure~\ref{hetero degen phase portrait figure}.  For the PDE \eqref{fully nonlinear elliptic pde}, the upstream and downstream states play the role of the equilibria in the original heteroclinic orbit, and the intermediate equilibrium would correspond to a distinct $x$-independent solution.  

To formulate this more precisely, observe that by composing with a sequence of translations in the $x$-direction, we can shift the incipient intermediate state upstream or downstream so that the solutions locally --- but not uniformly --- converge to a new front.   With that in mind, we say that a sequence of strictly monotone fronts $\{(u_n, \lambda_n)\}$ experiences a \emph{heteroclinic degeneracy} if there is a sequence $x_n \to \pm\infty$ so that the three limits
\begin{equation}    
  \lim_{x \to \mp\infty} \lim_{n \to \infty} u_n(x + x_n, \placeholder),
  \quad 
  \lim_{n \to \infty} \lim_{x \to +\infty} u_n(x, \placeholder),
  \quad 
  \lim_{n \to \infty} \lim_{x \to -\infty} u_n(x, \placeholder)  \quad  \textrm{exist and are distinct.} 
  \tag{A2} \label{gen hetero degen}
\end{equation}
One can further assume that $(u_n(\placeholder+x_n, \placeholder),\lambda_n)$ converges in $C_\loc^2(\overline{\Omega})$ to a monotone front $(u_*,\lambda_*)$.  

\begin{figure}
  \centering
  \includegraphics[scale=1,page=1]{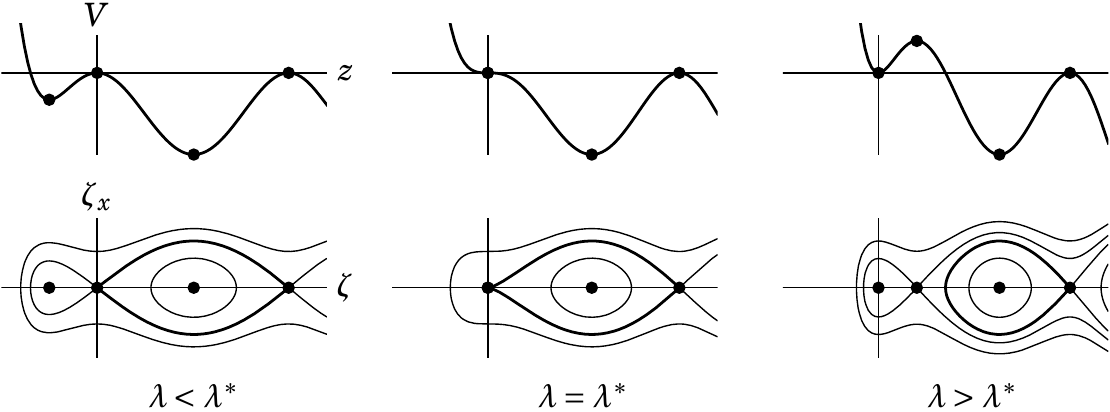}

  \caption{An ODE of the form \eqref{potential equation} that experiences a spectral degeneracy \eqref{gen ripples}. 
  }
  \label{ripples phase portrait figure} 
\end{figure}

Finally, it can happen that as we continue $\mathscr{K}_\loc$, the spectral non-degeneracy condition \eqref{ode spectral nondegen condition} fails.    Figure~\ref{ripples phase portrait figure} shows how this might occur in a specific example.   For $\lambda < \lambda^*$, there is a monotone increasing heteroclinic orbit connecting the constant solutions $Z_-=0$ and $Z_+=1$. This orbit persists for $\lambda = \lambda^*$, but the spectral condition \eqref{ode spectral nondegen condition} is violated at $Z_-$.  The Jacobian matrix for the corresponding planar system will then cease to be invertible downstream, and the orbit no longer decays exponentially as $x \to -\infty$.  For $\lambda > \lambda^*$, the heteroclinic orbit degenerates into a homoclinic orbit to $\zeta=1$, while $\zeta=0$ becomes a center.

The analogous scenario in the PDE setting should naturally involve  the spectrum of the linearized problem at infinity.  In particular, we say a sequence of strictly monotone fronts $\{(u_n,\lambda_n)\}$ experiences \emph{spectral degeneracy} if 
\begin{equation}
  \prineigenvalue^-(u_n, \lambda_n) \to 0 \qquad \textrm{or} \qquad \prineigenvalue^+(u_n, \lambda_n) \to 0.  \label{gen ripples} \tag{A3} 
\end{equation}
Recalling \eqref{spectral assumption},  we see that spectral degeneracy indicates resonance:   the essential spectrum of the linearized problem upstream or downstream moves through the origin.  Were this to occur, $\F$ will lose semi-Fredholmness and its zero-set may not be relatively compact.    In connection to traveling waves in reaction-diffusion equations, \eqref{gen ripples} corresponds to the onset of ``essential instability'' \cite{sandstede1999essential,sandstede2001essential}.  

\subsection{Statement of abstract results} \label{abstract results section}

Having developed the necessary intuition, we are now prepared to present the main global bifurcation theorem.

\begin{theorem}[Global bifurcation] \label{general global bifurcation theorem}
Consider the elliptic PDE \eqref{fully nonlinear elliptic pde}.  Let $\cm_\loc$ be a curve of strictly monotone front solutions which bifurcates from a singular point as in \eqref{local singular assumption} and satisfies the kernel \eqref{kernel assumption} and spectral \eqref{spectral assumption} conditions.  Then  $\mathcal{C}_\loc$ is contained in a global $C^0$ curve 
\[ \mathcal{C} := \left\{ \left(u(s), \lambda(s) \right) : 0 < s < \infty  \right\} \subset \F^{-1}(0)\]
of strictly monotone front solutions with the  properties enumerated below.
\begin{enumerate}[label=\rm(\alph*)]
\item \label{gen alternatives} \textup{(Alternatives)} For any sequence $s_n \to +\infty$, along some subsequence, $(u(s_n),\lambda(s_n))$, the blowup \eqref{gen blowup alternative}, heteroclinic degeneracy \eqref{gen hetero degen}, or spectral degeneracy \eqref{gen ripples} alternative will occur.  
\item \label{gen analyticity} \textup{(Analyticity)} At each point, $\iftcm$ admits a local real-analytic reparameterization.

  \item \label{gen reconnect} For all $s$ sufficiently large, $(u(s),\lambda(s)) \not\in \cm_\loc$.  In particular, $\cm$ is not a closed loop.
  \end{enumerate}
\end{theorem}

It bears repeating that the above theorem applies to a broad class of problems as it makes no structural hypotheses on the system beyond analyticity of $(A,B)$ and ellipticity.   There is a substantive body of work on fronts for semi-linear PDEs arising in reaction-diffusion equations (see \cite{volpert1994book} and the references therein).  To the best of our knowledge, however, Theorem~\ref{general global bifurcation theorem} is the first systematic treatment that applies even to fully nonlinear problems.  For example, in addition to the water wave applications discussed above, the general theory is used in a forthcoming paper to construct large nonlinear elastostatic fronts \cite{chen2020antiplane}.

Theorem~\ref{general global bifurcation theorem} is also distinctive in that it avoids making assumptions on the compactness properties of $\F$ beyond the local curve.  Classical global bifurcation theory makes comparatively stringent requirements that are appropriate for elliptic PDEs set on bounded domains but not the present problem.   For instance, Buffoni--Toland \cite{buffoni2003analytic} ask that the zero-set $\F^{-1}(0)$  be locally compact and $\F_u(u,\lambda)$ be Fredholm index $0$ for $(u,\lambda) \in \F^{-1}(0)$.  The seminal work of Rabinowitz \cite{rabinowitz1971some} assumes that $\F$ is locally proper and Fredholm index $0$ throughout its domain (though it need not be analytic). 

The basic philosophy inherent to our approach is that, on unbounded domains, it is more natural to think of the failure of these compactness properties as an \emph{alternative}, and then seek to classify it in terms of qualitative features of the solutions.  It is truly remarkable that the simple set of possibilities for the ODE \eqref{potential equation}, when properly interpreted, exhaustively categorize the limiting behavior for solution curves to the vastly more complicated PDE \eqref{fully nonlinear elliptic pde}.  

Let us conclude by briefly outlining how Theorem~\ref{general global bifurcation theorem} is used to construct large-amplitude bores.  Local curves $\cm_\loc$ of small-amplitude monotone front solutions to the internal wave problem \eqref{eqn:stream} were obtained in \cite{amick1989small,mielke1995homoclinic,makarenko1992bore,chen2019center}.  In \cite{chen2019center} this was done using a novel center manifold reduction method that is particularly well suited to verifying that the hypotheses \eqref{kernel assumption}, \eqref{spectral assumption}, and \eqref{local singular assumption} hold.  Due to its variational structure, \eqref{eqn:stream} possess several conserved quantities:  the mass flux, energy, and flow force through any vertical cross-section of the fluid domain must be the same.  The upstream and downstream states must therefore be \emph{conjugate} in that the values of these three quantities will agree.    For \eqref{eqn:stream}, this requirement is so restrictive that, in fact, at every $\lambda$, there is a unique downstream state that is conjugate to the fixed upstream state.  This insight drastically simplifies the task of computing the spectrum of the transversal linearized operators at $x =\pm\infty$, and indeed, we are able to rule out spectral degeneracy \eqref{gen ripples} entirely.  It also disqualifies the heteroclinic degeneracy alternative, as the three limiting states in \eqref{gen hetero degen} would be distinct and pairwise conjugate, which is impossible.    Thus blowup \eqref{gen blowup alternative} occurs as we follow the global bore curve.  Through elliptic regularity theory, we obtain uniform a priori bounds that show this must coincide with the stagnation limit \eqref{stagnation limit}.

\section*{Acknowledgments}

The research of RMC is supported in part by the NSF through DMS-1907584.  The research of SW is supported in part by the NSF through DMS-1812436.    A portion of this work was completed during a Research-in-Teams Program generously supported by the Erwin Schr\"odinger Institute for Mathematics and Physics, University of Vienna.

\bibliographystyle{siam}

\bibliography{projectdescription}

\def\cprime{$'$}
\begin{thebibliography}{10}

\bibitem{akers2014gravity}
{\sc B.~F. Akers, D.~M. Ambrose, and J.~D. Wright}, {\em Gravity perturbed
  {C}rapper waves}, Proc. R. Soc. Lond. Ser. A Math. Phys. Eng. Sci., 470
  (2014), pp.~20130526, 14.

\bibitem{ambrose2016global}
{\sc D.~M. Ambrose, W.~A. Strauss, and J.~D. Wright}, {\em Global bifurcation
  theory for periodic traveling interfacial gravity-capillary waves}, Ann.
  Inst. H. Poincar\'{e} Anal. Non Lin\'{e}aire, 33 (2016), pp.~1081--1101.

\bibitem{amick1987bounds}
{\sc C.~J. Amick}, {\em Bounds for water waves}, Arch. Rational Mech. Anal., 99
  (1987), pp.~91--114.

\bibitem{amick1989small}
{\sc C.~J. Amick and R.~E.~L. Turner}, {\em Small internal waves in two-fluid
  systems}, Arch. Ration. Mech. Anal., 108 (1989), pp.~111--139.

\bibitem{berestycki1992travelling}
{\sc H.~Berestycki and L.~Nirenberg}, {\em Travelling fronts in cylinders},
  Ann. Inst. H. Poincar\'{e} Anal. Non Lin\'{e}aire, 9 (1992), pp.~497--572.

\bibitem{buffoni2003analytic}
{\sc B.~Buffoni and J.~Toland}, {\em Analytic theory of global bifurcation},
  Princeton Series in Applied Mathematics, Princeton University Press,
  Princeton, NJ, 2003.
\newblock An introduction.

\bibitem{chen2019center}
{\sc R.~M. Chen, S.~Walsh, and M.~H. Wheeler}, {\em Center manifolds without a
  phase space for quasilinear problems in elasticity, biology, and
  hydrodynamics}, arXiv preprint arXiv:1907.04370,  (2019).

\bibitem{chen2020global}
\leavevmode\vrule height 2pt depth -1.6pt width 23pt, {\em Global bifurcation
  for monotone fronts of elliptic equations}, arXiv preprint arXiv:2005.00651,
  (2020).

\bibitem{chen2020antiplane}
\leavevmode\vrule height 2pt depth -1.6pt width 23pt, {\em Global bifurcation
  of anti-plane shear fronts}, Preprint,  (2020).

\bibitem{choi1996long}
{\sc W.~Choi and R.~Camassa}, {\em Long internal waves of finite amplitude},
  Physical Review Letters, 77 (1996), p.~1759.

\bibitem{choi1999fully}
\leavevmode\vrule height 2pt depth -1.6pt width 23pt, {\em Fully nonlinear
  internal waves in a two-fluid system}, J. Fluid Mech., 396 (1999), pp.~1--36.

\bibitem{constantin2016global}
{\sc A.~Constantin, W.~Strauss, and E.~V\u{a}rv\u{a}ruc\u{a}}, {\em Global
  bifurcation of steady gravity water waves with critical layers}, Acta Math.,
  217 (2016), pp.~195--262.

\bibitem{cordoba2016existence}
{\sc D.~C\'{o}rdoba, A.~Enciso, and N.~Grubic}, {\em On the existence of
  stationary splash singularities for the {E}uler equations}, Adv. Math., 288
  (2016), pp.~922--941.

\bibitem{crapper1957exact}
{\sc G.~D. Crapper}, {\em An exact solution for progressive capillary waves of
  arbitrary amplitude}, J. Fluid Mech., 2 (1957), pp.~532--540.

\bibitem{dias2003internal}
{\sc F.~Dias and J.-M. Vanden-Broeck}, {\em On internal fronts}, J. Fluid
  Mech., 479 (2003), pp.~145--154.

\bibitem{dyachenko2019folds}
{\sc S.~A. Dyachenko and V.~M. Hur}, {\em Stokes waves with constant vorticity:
  folds, gaps and fluid bubbles}, J. Fluid Mech., 878 (2019), pp.~502--521.

\bibitem{dyachenko2019stokes}
\leavevmode\vrule height 2pt depth -1.6pt width 23pt, {\em Stokes waves with
  constant vorticity: {I}. {N}umerical computation}, Stud. Appl. Math., 142
  (2019), pp.~162--189.

\bibitem{helfrich2006review}
{\sc K.~R. Helfrich and W.~K. Melville}, {\em Long nonlinear internal waves},
  in Annual review of fluid mechanics. {V}ol. 38, vol.~38 of Annu. Rev. Fluid
  Mech., Annual Reviews, Palo Alto, CA, 2006, pp.~395--425.

\bibitem{hur2020crapper}
{\sc V.~M. Hur and J.-M. Vanden-Broeck}, {\em A new application of {C}rapper's
  exact solution to waves in constant vorticity flows}, Eur. J. Mech. B Fluids,
  83 (2020), pp.~0--4.

\bibitem{hur2020exact}
{\sc V.~M. Hur and M.~H. Wheeler}, {\em Exact free surfaces in constant
  vorticity flows}, J. Fluid Mech., 896 (2020), p.~Paper No. R1.

\bibitem{kirchgassner1982wavesolutions}
{\sc K.~Kirchg{{\"a}}ssner}, {\em Wave-solutions of reversible systems and
  applications}, J. Differential Equations, 45 (1982), pp.~113--127.

\bibitem{Laget1997interfacial}
{\sc O.~Laget and F.~Dias}, {\em Numerical computation of capillary-gravity
  interfacial solitary waves}, J. Fluid Mech., 349 (1997), pp.~221--251.

\bibitem{makarenko1992bore}
{\sc N.~I. Makarenko}, {\em Smooth bore in a two-layer fluid}, in Free boundary
  problems in continuum mechanics ({N}ovosibirsk, 1991), vol.~106 of Internat.
  Ser. Numer. Math., Birkh{\"a}user, Basel, 1992, pp.~195--204.

\bibitem{maklakov2020note}
{\sc D.~V. Maklakov}, {\em A note on the existence of pure gravity waves at the
  interface of two fluids}, Phys. D, 401 (2020), pp.~132157, 5.

\bibitem{mielke1988reduction}
{\sc A.~Mielke}, {\em Reduction of quasilinear elliptic equations in
  cylindrical domains with applications}, Math. Methods Appl. Sci., 10 (1988),
  pp.~51--66.

\bibitem{mielke1995homoclinic}
\leavevmode\vrule height 2pt depth -1.6pt width 23pt, {\em Homoclinic and
  heteroclinic solutions in two-phase flow}, in Proceedings of the
  {IUTAM}/{ISIMM} {S}ymposium on {S}tructure and {D}ynamics of {N}onlinear
  {W}aves in {F}luids ({H}annover, 1994), vol.~7 of Adv. Ser. Nonlinear Dynam.,
  World Sci. Publ., River Edge, NJ, 1995, pp.~353--362.

\bibitem{miyata1988long}
{\sc M.~Miyata}, {\em Long internal waves of large amplitude}, in Nonlinear
  water waves, Springer, 1988, pp.~399--406.

\bibitem{pullin1988finite}
{\sc D.~I. Pullin and R.~H.~J. Grimshaw}, {\em Finite-amplitude solitary waves
  at the interface between two homogeneous fluids}, Phys. Fluids, 31 (1988),
  pp.~3550--3559.

\bibitem{rabinowitz1971some}
{\sc P.~H. Rabinowitz}, {\em Some global results for nonlinear eigenvalue
  problems}, J. Functional Analysis, 7 (1971), pp.~487--513.

\bibitem{sandstede1999essential}
{\sc B.~Sandstede and A.~Scheel}, {\em Essential instability of pulses and
  bifurcations to modulated travelling waves}, Proc. Roy. Soc. Edinburgh Sect.
  A, 129 (1999), pp.~1263--1290.

\bibitem{sandstede2001essential}
\leavevmode\vrule height 2pt depth -1.6pt width 23pt, {\em Essential
  instabilities of fronts: bifurcation, and bifurcation failure}, Dyn. Syst.,
  16 (2001), pp.~1--28.

\bibitem{spielvogel:overturning}
{\sc E.~R. Spielvogel}, {\em A variational principle for waves of infinite
  depth}, Arch. Rational Mech. Anal., 39 (1970), pp.~189--205.

\bibitem{sun2001existence}
{\sc S.~M. Sun}, {\em Existence of large amplitude periodic waves in two-fluid
  flows of infinite depth}, SIAM J. Math. Anal., 32 (2001), pp.~1014--1031.

\bibitem{toland1996stokes}
{\sc J.~F. Toland}, {\em Stokes waves}, Topol. Methods Nonlinear Anal., 7
  (1996), pp.~1--48.

\bibitem{turner1988broadening}
{\sc R.~E.~L. Turner and J.-M. Vanden-Broeck}, {\em Broadening of interfacial
  solitary waves}, Phys. Fluids, 31 (1988), pp.~2486--2490.

\bibitem{volpert1994book}
{\sc A.~I. Volpert, V.~A. Volpert, and V.~A. Volpert}, {\em Traveling wave
  solutions of parabolic systems}, vol.~140 of Translations of Mathematical
  Monographs, American Mathematical Society, Providence, RI, 1994.
\newblock Translated from the Russian manuscript by James F. Heyda.

\end{thebibliography}

\end{document}